\documentclass[article,letter,11pt]{article}
\usepackage[english]{babel}
\usepackage[T1]{fontenc}
\usepackage[bitstream-charter]{mathdesign}

\usepackage{kotex}
\usepackage{multicol}
\usepackage{amsmath}
\usepackage{amssymb}
\usepackage{amsthm}
\usepackage{graphicx}
\usepackage{xcolor}
\usepackage{authblk}
\usepackage{enumitem}
\usepackage{marvosym}
\usepackage{cite}
\usepackage{booktabs}
\usepackage{mathtools}
\usepackage{multicol}
\usepackage{multirow}
\usepackage[all]{nowidow}
\usepackage[font=small]{caption}	%footnotesize
\usepackage{tikz}
\usetikzlibrary{arrows.meta}
\usetikzlibrary{fadings}
\usepackage{smartdiagram}
\usepackage{metalogo}
\usepackage{vmargin}
\setmarginsrb{2cm}{2cm}{2cm}{2cm}{0mm}{0mm}{0mm}{10mm}
\usepackage{gensymb}
\usepackage[colorlinks=true, allcolors=blue, breaklinks=true]{hyperref}
\usepackage{apalike}
\theoremstyle{definition}

\usepackage{url}
\newcounter{example}[section]

\title{\vspace*{-1cm}\bfseries How can potential engineers seek and discover\\ their genuine vocations?}
\author{\normalsize Natanael Karjanto \thanks{\Letter: \texttt{natanael@skku.edu} \href{https://orcid.org/0000-0002-6859-447X}{\includegraphics[scale=0.08]{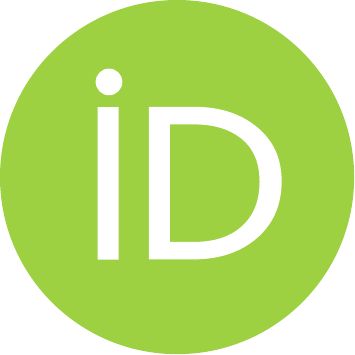}}}}
\affil{Department of Mathematics, University College, Natural Science Campus\\ Sungkyunkwan University, Suwon~16419, Republic of Korea}

\date{\vspace*{-0.5cm} \scriptsize Updated \today}

\begin{document}
\maketitle

\begin{abstract}
\noindent
The article describes what constitutes chemical engineering and its branches are in plain language. It outlines what to expect when one enrolls in a chemical engineering program as an undergraduate or graduate student. This may include core subjects to take, skillset to master, and other essential expertise that might be useful in the workplace. The discussion continues with career options for enthusiastic chemical engineers and how these young and early careers graduates could discover their reason for well-being and life purpose as aspiring chemical engineers.\\

\end{abstract}

% Section 1
\section{Introduction}

Engineering remains one of the favorite majors among many college applicants and engineering graduates continue to gain respect and dignity among professional employees, with handsome take-home salaries and high job satisfaction reports~\cite{gonzales2018is,perry2020job}. According to the \emph{City \& Guilds Happiness Index} survey from the UK, chartered engineers are the happiest professional around, they beat lawyers, scientists, architects, and accountants when it comes to job satisfaction~\cite{anscombe2005they}. What is particularly interesting, according to a survey by \emph{Payscale.com} in 2010, the chemical engineering major was ranked number one in terms of job satisfaction, while psychology was ranked last among majors in the survey~\cite{light2010psych}.

The field of engineering itself is certainly a broad discipline, and we will focus only on one of the main branches of engineering, i.e., chemical engineering. We seek to answer the following research questions:
\begin{itemize}[leftmargin=1em]
\item What is chemical engineering? And what does it cover?
\item What are the core subjects that one needs to master if he/she wants to become a chemical engineer?
\item What kind of careers that newly trained chemical engineers could pursue?
\end{itemize}
In what follows, the remainder of this introduction will answer the first research question, review the body of literature, and lay out a theoretical framework for our discussion.

Chemical engineering is a certain type of and one of the main branches of engineering that deals with the study of operation and design of chemical plants as well as methods of improving production. It applies practices in mathematics and scientific principles for the design, development, and evaluation of operational systems. This includes employing the application of physics, chemistry, biology, and engineering principles harmoniously to successfully perform its operations. Chemical engineering contributes to the development of a wide range of new as well as improved products and processes. These comprise strong materials to resist extreme temperatures, new fuels for reactors, and booster propulsion, medicines, vaccines, serum, and plasma for mass distribution, among others~\cite{beakley1979careers}.

Chemical engineering is both an art and science, in addition to engineering itself certainly. Chemical engineering as science means that some parts of the field are thoroughly understood theoretically. As an art, there are other areas of chemical engineering that are only partially understood theoretically. Similar to other engineering fields, chemical engineering is an art to a certain extent since often, chemical engineers conduct their work based on experience and judgment~\cite{badger1955introduction}.

The profession of chemical engineering deals with the technology of the chemical and process industries. Chemical engineers must be able to develop, design, and engineer both the equipment and process involved in their vocation. In addition to being able to operate plants efficiently, safely, and economically, chemical engineers must also understand the functions and characteristics of the product outputs~\cite{badger1955introduction}.

Chemical engineers carry out chemical processes on a commercial scale by converting raw materials into useful products. They must be able to apply scientific and engineering principles when designing and operating plants for material production. Chemical engineers are also involved in various aspects of plant design and operation, process design and analysis, chemical reaction engineering, materials synthesis and processing, environmental aspects, as well as safety and hazard assessments. The areas that chemical engineers need to deal with encompassing the manufacture of commodity chemicals, specialty chemicals, petroleum refining, microfabrication, fermentation, synthetics and plastics, and biomolecule production. 

The body of literature contains rich and diverse coverage on chemical engineering. Although what we discuss in this article is far from exhaustive, they provide a basis for understanding chemical engineering nonetheless. Interested readers are encouraged to consult the bibliography listed at the end of this article and consult the references therein. 

Kim (2002) covered chemical engineering from a historical point of view~\cite{kim2002chemical}. Hipple discussed basic concepts of chemical engineering in an easy-to-understand way, which enables non-chemical engineers to understand fundamental concepts of chemical processing, design, and operation effortlessly~\cite{hipple2017chemical}. Nnaji (2019) offered a comprehensive overview of the evolution, essence, concept, principles, functions, and applications of chemical engineering~\cite{nnaji2019introduction}. Ogawa (2007) attempted to enlighten the usefulness of information entropy for many phenomena that appear in chemical engineering~\cite{ogawa2017chemical}. 

In particular, for high school students who are considering majoring in chemical engineering, Stimus (2013) provided a very simple and non-technical introduction to chemical engineering~\cite{stimus2013beginner}. The book was particularly written for readers who have no background and experience in the field. Additionally, also written with high school students and recent graduates in mind, Ridder (2016) arranged a gentle tour through what chemical engineers conduct in practice, how they accomplish it, what kinds of employments and careers are open to them, and how much salary they can expect~\cite{ridder2016balancing}.
\tikzfading[name=fade inside,
inner color=transparent!10,
outer color=transparent!40] 
\begin{figure}[h]
\begin{center}
\begin{tikzpicture}
\shade[ball color=green!80!blue,path fading=fade inside] (0,0) circle (1.8);
\begin{scope}[xshift=-2cm,transform canvas={rotate=-180}]
\shade[ball color=blue!80!green,path fading=fade inside] (0,0) circle (1.8);
\end{scope}
\begin{scope}[xshift=1cm,yshift=2cm]
\shade[ball color=red!80!orange,path fading=fade inside] (0,0) circle (1.8);
\end{scope}
\node at (30:5.5) {Chemical engineering};
\node at (180:3)  {Education};
\node at (0:5) 	  {Psychology};
\end{tikzpicture}
\end{center}
\caption{A theoretical framework for meaningful and satisfying careers for aspiring chemical engineerss.}	\label{framework}
\end{figure}

The focus of this article is on both the educational and psychological aspects of chemical engineering. The theoretical framework for our discussion is an intersection of chemical engineering with education and psychology, as shown in Figure~\ref{framework}. After this introduction, the article is organized as follows. Section~\ref{edu} covers the curriculum and educational aspects of chemical engineering. Section~\ref{car} discusses possible career options for aspiring chemical engineers. Finally, Section~\ref{con} concludes our discussion.

% Section 2
\section{Chemical engineering education}	\label{edu}

\subsection{General curriculum}

According to Thompson and Ceckler (1977), students who take their first course in an introduction to chemical engineering should understand the nature and scope of the chemical process industry, comprehend physical and chemical principles as they relate to the analysis of the chemical processes, and proficient in a systematic approach to organizing information and solving problems~\cite{thompson1977introduction}.
\begin{table}[h]
{\small 
\caption{An example of curriculum for an undergraduate program in chemical engineering. Course names are adapted from the University of Virginia (Charlottesville, Virginia), Clemson University (Clemson, South Carolina), University of Colorado (Boulder, Colorado), and Auburn University (Auburn, Alabama).}		\label{curriculum}
\begin{center}
\vspace*{-10pt}	
\begin{tabular}{@{}lcclc@{}}
\toprule	
\multicolumn{5}{c}{Freshman} \\ \hline %\midrule
Fall 						& Credits & {\qquad} 	& Spring 				& Credits 	\\ \cline{1-2} \cline{4-5}
Calculus 1					& 4		  &				& Calculus 2 			& 4			\\
General Physics 1     		& 4		  &				& General Physics 2 	& 4			\\
General Chemistry 1   		& 4       &				& General Chemistry 2	& 4			\\
Introduction to Engineering & 3		  &				& Biology for Engineers	& 3			\\
HSS Elective				& 3		  &				& HSS Elective 			& 3			\\ \cline{2-2} \cline{5-5}
                            & 18      &             &                       & 18        \\ \midrule 
\multicolumn{5}{c}{Sophomore} \\ \hline %\midrule                             
Fall 						& Credits & {\qquad} 	& Spring 				& Credits 	\\ \cline{1-2} \cline{4-5}
\multirow{2}{*}{Calculus 3} & \multirow{2}{*}{4} &	& Differential Equations  			& \multirow{2}{*}{4}			\\
							&		  &				& and Linear Algebra    & \\	
Principles of Chemical Engineering & 4&				& Fluid Mechanics 		& 3			\\
Organic Chemistry 1   		& 4       &				& Organic Chemistry 2	& 4			\\
Materials and Energy Balance& 3		  &				& Physical Chemistry	& 3			\\
HSS Elective				& 3		  &				& Thermodynamics		& 3			\\ \cline{2-2} \cline{5-5}
							& 18      &             &                       & 17        \\ \midrule	
\multicolumn{5}{c}{Junior} \\ \hline %\midrule                             
Fall 						& Credits & {\qquad} 	& Spring 				& Credits 	\\ \cline{1-2} \cline{4-5}
Heat Transfer  				& 3 	  &				& Mass Transfer 		& 3			\\
Physical Chemistry 			& 3		  &				& Chemical Reaction		& 3			\\
Separation Process   		& 3       &				& Chemical Engineering Analysis	& 3			\\
Statistical Methods			& 3		  &				& Engineering Materials	& 3			\\
Unrestricted Elective		& 3		  &				& Kinetics				& 3			\\ \cline{2-2} \cline{5-5}
							& 15      &             &                       & 15        \\ \midrule		
\multicolumn{5}{c}{Senior} \\ \hline %\midrule                             
Fall 						& Credits & {\qquad} 	& Spring 				& Credits 	\\ \cline{1-2} \cline{4-5}
Process Control				& 3 	  &				& Process Design 		& 3			\\
Physical Chemistry 			& 3		  &				& Chemical Reaction		& 3			\\
Process Safety		   		& 3       &				& Bioprocess Engineering& 3			\\
Major Elective 				& 3		  &				& Major Elective		& 3			\\
Unrestricted Elective		& 3		  &				& Unrestricted Elective	& 3			\\ \cline{2-2} \cline{5-5}
							& 15      &             &                       & 15        \\ \midrule		
\bottomrule
\end{tabular}
\end{center}
}
\end{table}

In order to become a chemical engineer graduate, one needs to obtain a competent credential. Pursuing a bachelor's degree in chemical engineering is the common path to follow, which usually takes four years to complete. One may continue to graduate schools to obtain a master's and/or PhD degree to delve further. The former commonly requires two years of study while the latter may range from three to six years to complete, depending on various factors. In addition to acquiring important skills, adequate experiences, and necessary credentials, these degrees open the door to a wide array of exciting academic and careers opportunities, as we will discuss in Section~\ref{car}. 

In addition to enrolling in required core courses, students enrolled in an undergraduate program of chemical engineering also need to take general courses in mathematics and sciences, which may include sequences of calculus, physics, chemistry, and biology courses. The calculus courses cover differential, integral, multivariable, vector calculus, linear algebra, and differential equations. The physics modules include classical mechanics, thermodynamics, electricity, magnetism, and special relativity. It should not be surprising that in addition to two-semesters general chemistry courses covering atomic and molecular theory, chemical periodicity, stoichiometry, chemical reactions, chemical kinetics, and chemical equilibrium, students also need to take a sequence of courses in organic chemistry. Some institutions also require their students to complete a set of social science and humanities elective courses.

Common themes covered in chemical engineering programs include physics, chemistry, kinetics, chemical reactions, separation processes, electrochemical systems, fluid mechanics, transport phenomena, energy conservation, and many others. Many of these subjects are often accompanied by laboratory activities and computer applications. Table~\ref{curriculum} displays an example of a curriculum for an undergraduate program in chemical engineering.

\begin{figure}[h]
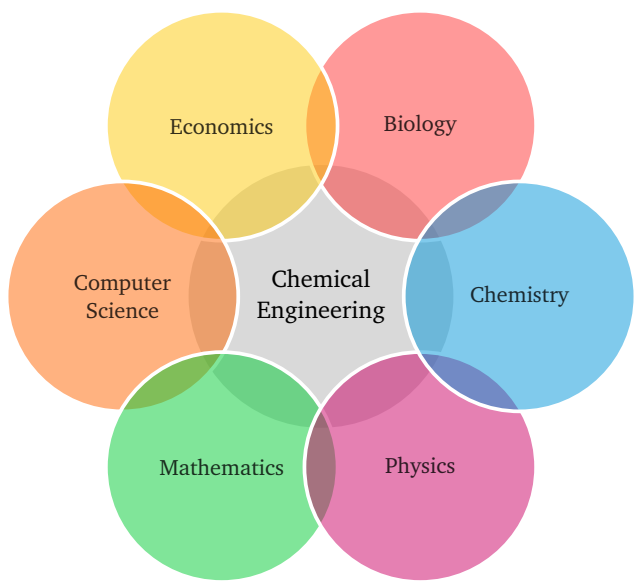

\begin{center}
\resizebox{.5\linewidth}{!}{%
\smartdiagramset{
	bubble node size =3.5cm, 
	bubble center node font = \normalsize,
	bubble node font = \small, 
	distance center/other bubbles = 1cm, 
	set color list = {red!80!white, yellow!80!red, orange!80!red, green!80!blue, purple!80!magenta, cyan!80!blue}
}%
\smartdiagram[bubble diagram]{%
	Chemical\\ Engineering,
	Biology,
	Economics,
	Computer\\Science,
	Mathematics,
	Physics,
	Chemistry%
}%
}% end resizebox
\end{center}
\caption{Chemical engineering consists of applications of mathematics, physics, chemistry, biology, computer science, and economics.}	\label{intersection}
\end{figure}

It is a common misconception among prospective undergraduate and high school students that they can learn more about and specialize in chemistry once they take a chemical engineering program. While the notion is not entirely wrong, it turns out that only around 20\% of the courses are related to chemistry, particularly physical and organic chemistry, as we observe in Table~\ref{curriculum}. Similar to other engineering degrees, mathematics and physics are ultimately the most important elements of the program. There are also some courses related to biology, but these would typically be for more specialized tracks. See Figure~\ref{intersection}.

To see the difference between chemistry and chemical engineering, it is worth understanding the difference between science and engineering itself. Although they are different, both are closely related and have existed since the dawn of history. Science involves expanding a body of knowledge about a particular topic by conducting a sequence of simulations or experiments. In contrast, engineering applies those scientific ideas to real life. While scientists use scientific methods in their research, engineers employ scientific principles to design equipment, construct tools, and realize products~\cite{stimus2013beginner}.

By examining the curricula closely, there are significant differences between chemistry and chemical engineering programs. Studying and majoring in chemistry means exploring large repositories of chemical reactions and investigating the fundamental theories behind their existence. In addition to learning in great detail the theory behind numerous instruments they use for chemical analysis, chemistry majors also study and experiment a lot with chemical synthesis, i.e., construction of complex chemical compounds from simpler ones. On the other hand, majoring in chemical engineering focuses on the practical implementation of various reaction patterns and relates how different transport phenomena affect the outcome of those reactions~\cite{ridder2016balancing}.

\subsection{Mathematics and physics in chemical engineering}

When it comes to mathematics and physics contents in chemical engineering, the amount is illuminating. In addition to learning a sequence of calculus courses, linear algebra, differential equations, and statistical methods, prospective chemical engineers also study mathematically-intensive physics modules, including chemical thermodynamics, fluid mechanics, mass transfer, transport phenomena, and reaction engineering. In this subsection, after a brief literature coverage, we provide one classic example in fluid mechanics where mathematics and physics come alive and discover its applications in chemical engineering, i.e., the classical Graetz problem~\cite{graetz1882ueber,graetz1885ueber}.

As always, the body of literature will never be exhaustive, and we are barely scratching its surface in this article. Mickley et al. (1957) has written a classic on applied mathematics in chemical engineering, and it is regarded as an excellent textbook for chemical engineers who are willing to perform their modeling~\cite{mickley1957applied}. Rice and Do (2012) demonstrated how classical mathematics solves a broad range of new application problems in chemical engineering in their book's second edition~\cite{rice2012applied}. Focusing on mathematics application in chemical engineering, Loney (2016) also  addressed the setup and verification of mathematical models using experimental data in his expanded and revised version monograph~\cite{loney2016applied}.

Acrivos and Admundson (1955) demonstrated a brief survey on applications of matrix algebra and calculus in chemical engineering problems~\cite{acrivos1955applications}. Sweeny et al. (1964) reviewed mathematics applications in chemical engineering, particularly in statistics, basic and industrial operations research, and optimization~\cite{sweeny1964mathematics}. Ramkrishna and Amundson (2004) provided an overview of mathematics in chemical engineering during the second half of the past century~\cite{ramkrisna2004mathematics}. Wang et al. (2013) observed a common mathematical feature of geometric and parametric singularities observed in their research in chemical engineering~\cite{wang2013celebrating}. Finlayson et al. (2015) comprehensively covered mathematical topics useful for chemical engineering~\cite{finlayson2015mathematics}. 
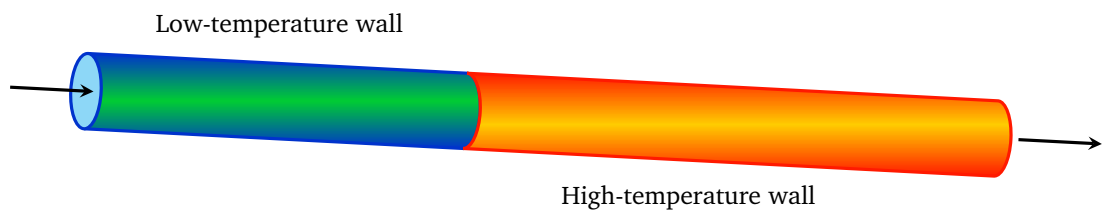
\begin{figure}[h]
\begin{center}
\begin{tikzpicture}[rotate=-3]
\node at (2.5,1) {\small Low-temperature wall};
\node at (8,-1) {\small High-temperature wall};
\draw [blue!80!green, very thick, top color=blue!80!green,
bottom color=blue!80!green,
middle color=green!80!blue] (0,-0.5) -- ++(5,0) 
arc(-90:90:0.2 and 0.5) -- ++(-5,0) ;
arc(90:-90:0.2 and 0.5)--cycle;
\draw[blue!80!green, very thick, fill=cyan!40] (0,0) ellipse(0.2 and 0.5);
\draw [red!80!orange, very thick, top color=red!80!orange,
bottom color=red!80!orange,
middle color=yellow!80!red] (5,-0.5) -- ++(7,0) 
arc(-90:90:0.2 and 0.5) -- ++(-7,0) 
arc(90:-90:0.2 and 0.5)--cycle;
%\draw[blue!80!green, very thick, fill=cyan!10] (0,0) ellipse(0.2 and 0.5);
\draw[->,>=stealth,very thick] (-1,0) -- + (1.1,0);
\draw[->,>=stealth,very thick] (12.3,0) -- (13.4,0);
\end{tikzpicture}
\end{center}
\caption{A schematic diagram for the classical Graetz problem. Fluid enters upstream from the left-hand side of the tube encounters a low-temperature wall. As it travels downstream the tube, it experiences a sudden jump in heat as the surrounding tube possesses a higher temperature surface.}		\label{graetz}
\end{figure}

The classical Graetz problem is a fundamental tube flow problem that couples fluid flow with heat and/or mass transfer. It considers the thermal entry of an incompressible fluid in a circular pipe with a fixed velocity profile. There exists a sudden jump in the tube temperature, with lower and higher temperatures at the downstream and upstream positions, respectively. The fluid enters the tube with a particular temperature at the upstream position, which can be higher or lower than the wall temperature. See Figure~\ref{graetz}. The objective is to determine the behavior of the temperature profile as it travels along the tube downstream. Since the flow is incompressible, the velocity distribution does not depend on the varying temperatures. In particular, for laminar flow in a circular pipe, the velocity profile is parabolic~\cite{shah1978laminar,kee2003chemically}. The Graetz problem finds applications in chemical reactors, heat exchangers, blood flow, and other fluid transport phenomena.

By calling the boundary-layer assumptions for neglecting the axial conduction and diffusion, and by employing the parabolic velocity profile, the governing equation for the temperature can be expressed in the following partial differential equation (PDE) in nondimensional quantities:
\begin{equation*}
\frac{\partial T}{\partial z} = \frac{1}{(1 - r^2)} \, \frac{1}{r} \, \frac{\partial}{\partial r} \left(r \frac{\partial T}{\partial r} \right),
\end{equation*}
with the following boundary conditions:
\begin{equation*}
T(r = 0, z) = 1, \qquad T(r = 1, z) = 0, \qquad \text{and} \qquad \frac{\partial T}{\partial r}(r = 0,z) = 0.
\end{equation*}
By implementing the method of separation of variables, i.e., writing $T(r,z) = R(r)Z(z)$, the PDE is decomposed into two ordinary differential equations (ODEs). In particular, $R$ satisfies the following Sturm-Liouville boundary value problem (BVP):
\begin{align*}
-\frac{1}{r} \frac{d}{dr} \left(r \frac{dR}{dr} \right) &= \lambda^2 (1 - r^2) R \\
R(1) = 0, \qquad &\text{and} \qquad \frac{dR}{dr}(0) = 0.
\end{align*}
Performing another set of transformations to both dependent and independent variables, i.e., $x = -\lambda r^2$ and $R = e^{-x/2} y(x)$, we obtain the well-known Kummer's differential equation~\cite{huang1984heat}:
\begin{equation*}
x \frac{d^2y}{dx^2} + (1 - x) \frac{dy}{dx} - \left(\frac{1}{2} + \frac{\lambda}{4} \right) y = 0.
\end{equation*}
This Kummer's ODE belongs to a class of confluent hypergeometric equation, a degenerate form of hypergeometric ODE where two of the three regular singularities merge into an irregular singularity~\cite{kumer1837de,slater1960confluent}. Although a second-order ODE yields two linearly independent solutions, it is only one that is bounded at $x = 0$. By requiring that both $R(0)$ and $y(0)$ must remain finite, we rule out the singular solution, and hence, we acquire only the first kind of Kummer's function. Up to a multiple of constant, this regular solution can be expressed in the following form: 
\begin{equation*}
y(x) = M\left(\frac{1}{2} - \frac{\lambda}{4}, 1, x \right) = {\,}_1 F_{1}\left(\frac{1}{2} - \frac{\lambda}{4}; 1; x \right)
= \sum_{n = 0}^{\infty} \left(\frac{1}{2} - \frac{\lambda}{4}\right)^{(n)} \frac{x^n}{1^{(n)} n!},
\end{equation*}
where the Pochhammer function $\xi^{(n)}$, $n \in \mathbb{N}_0$, (also known as the rising or ascending factorial) is defined as
\begin{align*}
\xi^{(0)} &= 1, \\
\xi^{(n)} &= \prod_{k = 0}^{n - 1} (\xi + k) = \xi (\xi + 1) (\xi + 2) \cdots (\xi + n - 1).
\end{align*}
By applying the boundary condition at the tube wall, i.e., $R(1) = 0$ or $y(1) = 0$, we obtain the following transcendental equation for the eigenvalues, or rather, the square root of eigenvalues:
\begin{equation*}
 M\left(\frac{1}{2} - \frac{\lambda_n}{4}, 1, -\lambda_n \right) = 0, \qquad n \in \mathbb{N}_0.
\end{equation*}
For each value $\lambda_n$, there exists an associated eigenfunction $R_n(r)$ given by
\begin{equation*}
R_n(r) = e^{-\frac{1}{2}\lambda_n r^2} y_n(-\lambda_n r^2).
\end{equation*}
A general solution to the classical Graetz problem can then be obtained by applying the linear superposition principle. Furthermore, the unknown coefficients can be retrieved by utilizing the orthogonality properties of the Sturm-Liouville problem after employing the initial condition accordingly~\cite{huang1984heat}. The readers who are interested in delving more about the Graetz problem may consult~\cite{notter1972solution,michelsen1974graetz,papoutsakis1980extended,papoutsakis1980dirichlet,baron1997graetz}; for confluent geometric Kummer's function and equation, review~\cite{slater1960confluent,buchholz1969confluent,georgiev2003new,georgiev2005kummer,ancarani2008derivatives}. For more information on the Sturm-Liouville problem, see the listed monographs and their references therein~\cite{pryce1993numerical,amrein2005sturm,gwaiz2008sturm,zettl2010sturm,teschl2021ordinary,haberman2013applied,kravchenko2020direct,zettl2021recent}.

% Section 3
\section{Careers in chemical engineering}	\label{car}

\subsection{Career and employment}

A degree in chemical engineering is a path in making a positive difference in this world. There are numerous facts and inspiring stories about how chemical engineers have shaped our world over the years. Today chemical engineers have benefited from systematic subjects and satisfying careers thanks to the early pioneers in the field and thus they are standing on the shoulders of giants, so to speak.

For example, the British chemist George Edward Davis (1850--1907) is regarded as the founding father of chemical engineering. In addition to identifying extensive characteristics commonly found in the chemical industry, Davis published lecture series where he defined chemical engineering as a discipline and wrote an influential handbook in chemical engineering~\cite{davis1904handbook,freshwater1980george,freshwater2004george,flavel2012meet}. An American chemist Arthur Dehon Little (1863--1935) introduced and developed the concept of unit operations used in physical and chemical changes and promoted the set of activities of industrial research and technological development~\cite{keyes1937arthur,servos1980industrial,servos1996physical}. The first woman to earn a PhD in chemical engineering and the first female member of the American Institute of Chemical Engineers, Margaret Hutchinson Rousseau (1910--2000) designed the first commercial penicillin production plant~\cite{cep2008fifty,hatch2006changing,madhavan2015think}. Although polyethylene was discovered twice, it was a chemical engineer from England, Dermot Manning, who built high-pressure research reactors that made experiments and full-scale production possible~\cite{flavel2011dermot,michels2018how,clift2019managing}.

While not every chemical engineering graduate will make a significant breakthrough like these people, many graduates still have a lot of opportunities to work and contribute in a wide range of fields and industries. Their roles involve not only developing and improving existing processes but also discovering and creating new techniques for altering materials and producing useful products. Other equally meaningful sectors which utilize the skills developed in a chemical engineering degree include quality assurance, consultancy, and manufacturing. A few popular examples of industries where chemical engineers work and progress in their careers are petroleum and natural gas, pharmaceuticals, processed food, pulp and paper, metallurgy, mining, and metal extraction, household products, polymers, and medical devices~\cite{ridder2016balancing}.

The potential chemical engineers and many readers might be pleased to understand that skillful graduates with a degree in chemical engineering earn handsome salaries and enjoy generous benefits. If money is not the main concern, a lofty salary may also mean that chemical engineers deserve the respect and dignity for their qualified skillsets and credentials~\cite{sandel2020tyranny}. In the United States, the starting salaries for chemical engineers rank among the highest, approaching \$70,000 per annum, just below petroleum and computer engineers' salaries~\cite{clark2013chemical}. A more recent report from the US Bureau of Labor Statistics suggests that the median annual salary for chemical engineers stood at \$108,540~\cite{bls2020outlook}. Depending on industries, the wage may vary significantly, with the top one hoards more than \$150,000 annually. The top-paying industries for chemical engineers are oil and gas extraction, management of companies and enterprises, petroleum and coal products manufacturing, computer and peripheral equipment manufacturing, and local government, albeit excluding schools and hospitals~\cite{bls2020gov}. 

\subsection{Vocation and \emph{ikigai}}

While working in a particular industry might attract particularly high salaries, and nothing wrong with that, many new and young chemical engineer graduates might also want to seek their meaning and purpose in life while they pursue a career in chemical engineering. It is an indisputable fact that although money is not everything, many things require money. Despite the conventional wisdom that money cannot buy happiness, well-known research revealed that wealthier people tend to experience life satisfaction and feel happier than the unfortunate ones~\cite{kahneman2010high,jebb2018happiness,donnelly2018amount,sirgy2021dual}. Moreover, another set of studies demonstrates that richer people also tend to be healthier and live longer than their poorer counterparts~\cite{pritchett1996wealthier,bosworth2016later,zaninotto2020socio,finegood2021association}.

In his best-selling book, Peter Buffet, the youngest son of the philanthropist and billionaire investor Warren Buffet, talked about making a life by identifying and pursuing passions, determining goals, and defining success. Since life is what we make it, nobody else can and should do it for us. Although we might think that Peter has enjoyed a life of endless privilege, he admitted that the only real inheritance handed down from his wealthy parents was neither property nor saving account, but a philosophy instead: ``Find your bliss. Forge your path in life.'' Whilst choosing an independent path as a musician and producer, Peter not only followed his passions and established his identity, but also reaped his successes as well~\cite{buffett2011life}.

In particular, he discussed at length about \emph{vocation}. Vocation can be defined as ``life's work'', ``a trade or profession'', ``a person's employment or main occupation, especially regarded as particularly worthy and requiring great dedication'', ``a strong feeling of suitability for a particular career or occupation''. The word comes from the Latin word \emph{vocatio} or \emph{vocare}, which means ``to call'' or ``to summon''. Both terms \emph{calling} and \emph{vocation} are often applied to a sense of purpose or direction that leads an individual toward some kind of personally fulfilling and/or socially significant engagement within the work role. Despite its commonly used in non-religious contexts, the meaning of these expressions has roots in Christianity~\cite{schuurman2004vocation,dik2009calling,muller2017dictionary}. 

Buffet Junior offered two definitions of vocation, one to mean ``a passion for the work we do'', another encompasses a broader perspective, ``the tug we feel toward the life that is right for us, the life that is truly our own.'' He also admitted that in addition to being mysterious, life vocations are hardly linear. Embracing a vocation can be so scary, particularly for inexperienced ones. Recognizing and welcoming a calling raises the stakes in life, this is the reason why it can be so tough for young people to concede and commit to their true vocation. While hunches and inklings might spark along the way, we need to be patient and open-minded in identifying specific ways to conjugate them with our particular temperament and abilities to upgrade them to become a true vocation. Finally, Peter also solved the philosophical dichotomy between being versus doing. For an individual with a genuine work-vocation, doing equals being; whereas for a person with a sincere happiness-vocation, being equals doing~\cite{buffett2011life}. 

When discovering meaningful goals and finding true purpose in life, we can learn a lot from our Japanese counterparts. Japan has one of the highest life expectancies in the world, with an average of 84.79 years in 2021, a record high of 87.74 and 81.64 years for women and men, respectively. Although many factors are involved~\cite{buettner2012blue,tsugane2021why}, there is a unique concept in particular, i.e., \emph{ikigai}. It is a Japanese mindset of well-being referring to something that gives a person a sense of purpose, a reason for living and aspiration in life. They even talk about it as a reason to get up in the morning. Whilst referring to Japanese culture, a dictionary defines ikigai as `` motivating force; something or someone that gives a person a sense of purpose or a reason for living'', a general definition can mean ``something that brings pleasure or fulfillment''~\cite{oed2021}. 

The word is ikigai (生き甲斐) literally means ``a reason for being alive''. The term compounds two Japanese words: \emph{iki} (生き), meaning ``to live, alive, life'' and \emph{kai} (甲斐), meaning ``worth, avail, effect, result, reason, benefit''. The literal meaning is comparable to the borrowed French phrase \emph{raison d'être} (reason to be). In addition to describing having a sense of purpose in life and being motivated, ikigai also means the feeling of accomplishment and fulfillment when people pursue their passions~\cite{mathews1996stuff,mogi2017little,garcia2017ikigai,schippers2017ikigai,kumano2018concept,schippers2019life}. 

After reading this article, we hope that our readers and potential chemical engineers will have something to contemplate as they decide what they are going to do with their life. We finish this section with powerful quotes from a published author Robert (Bob) A. Moawad (1941--2007) and a writer and humorist Mark Twain (1835--1910). Both are Americans. The former said, ``The best day of your life is the one on which you decide your life is your own. No apologies or excuses. No one to lean on, rely on, or blame. The gift is yours--it is an amazing journey--and you alone are responsible for the quality of it. This is the day your life really begins.''  The latter thought-provokingly stated, ``The two most important days in your life are the day you are born and the day you find out why.''

% Section 4
\section{Conclusion}						\label{con}

We have introduced chemical engineering to general readers in this article. Our discussion encompasses the definition and scope of chemical engineering, its general curriculum at the undergraduate level, and potential career options for aspiring chemical engineers. For each topic, we have provided selective yet important literature to support the argument. We hope this article will illuminate not only many high school students who contemplate applying to an engineering program, but also parents, teachers, and academic counselors who desire their children, pupils, and advisees to uncover their great but might be hidden potentials.

As we have observed from the educational aspect, the curriculum of a chemical engineering undergraduate program is a balanced combination between mathematics, physics, and chemistry. Hence, anyone who is good at them and is interested in delving further into the application aspects of these subjects might want to consider majoring in chemical engineering. A classic example of the Graetz problem that marries fluid flow and heat/mass transfer demonstrates this feature. Those who are interested in specializing further, possibly at an advanced undergraduate or graduate level may discover multidisciplinary aspects of the subject, with biology, computer science, and economics among others appear to become significantly important.

Graduates of chemical engineering work in various industries and sectors, including production and service, public and private. Chemical engineers solve diverse and challenging yet interesting problems under various technical, economic, and environmental constraints. They earn respect and dignity thanks to their professional skillsets and credentials, as shown by attractive remuneration offers and packages. With future global challenges in healthcare, sustainability, and security, current and future chemical engineers have plenty of opportunities to contribute to the world for a better society.

% References
\end{document}